\documentclass[11pt]{article}
\usepackage{geometry}                
\usepackage[parfill]{parskip}    
\usepackage{graphicx}
\usepackage{amssymb}
\usepackage{epstopdf}
\usepackage{xcolor}
\usepackage{amsthm}
\usepackage{amsmath}
\usepackage{comment}
\usepackage{hyperref}
\usepackage{tikz}
\newenvironment{proof1}[1]{\begin{trivlist} \item[] {\em Proof of #1:}}{\newline \textcolor{white}{.}\hfill $\Box$
                      \end{trivlist}}
\theoremstyle{plain}

\newtheorem{thm}{Theorem}[section]
\newtheorem{lem}{Lemma}[section]
\newtheorem{defn}{Definition}[section]

\newtheorem{cor}{Corollary}[section]

\newtheorem{rem}{Remark}

\newcommand{\pa}{\partial}

\newcommand{\eps}{\epsilon}

\vspace{-0.3in}
\title{An isoperimetric inequality for surfaces formed from spherical polygons}
\date{\today}\author{Farhan Azad, Thomas Beck, and Karolina Lokaj}
\begin{document}
\maketitle

\begin{abstract}

 We {give a new proof of} an isoperimetric inequality for a family of closed surfaces, which have Gaussian curvature identically equal to one wherever the surface is smooth. These surfaces are formed from a convex, spherical polygon, with each vertex of the polygon leading to a non-smooth point on the surface. For example, the surface formed from a spherical lune is a surface of revolution, with two non-smooth tips. {Combined with a straightforward approximation argument, this inequality was first proved by B\'erard, Besson, and Gallot in \cite{BBG}, where they provide a generalization of the L\'evy-Gromov isoperimetric inequality.} The inequality implies an isoperimetric inequality for geodesically convex subsets of the sphere, and, using a Faber-Krahn theorem, it also implies a lower bound on the first Dirichlet eigenvalue of a region of a given area on the closed surfaces. Via approximation, we convert this into a lower bound on the first Dirichlet-Neumann eigenvalue of domains contained in geodesically convex subsets of the sphere. 

\end{abstract}

\section{Introduction}

An isoperimetric inequality on a Riemannian surface $S$  provides a lower bound on the length of a closed curve enclosing a given surface area. Denoting $L$ to be the length of the curve and $A$ to be the enclosed surface area, a classical example of such an inequality is
\begin{align*}
    L^2 \geq 4\pi A - K A^2.
\end{align*}
This inequality holds when $S$ is a simply-connected surface of constant sectional curvature $K$, and so it holds for the Euclidean plane ($K=0$), the sphere ($K>0$), and the hyperbolic plane ($K<0$), and was first proved in this generality in \cite{Sc1}. For all $K\in\mathbb{R}$,  equality is obtained precisely when the region is a geodesic disc, and so these are the unique perimeter-minimizing regions. 

In the variable curvature case, Benjamini and Cao \cite{BC} have proved an isoperimetric inequality for a class of rotationally symmetric curved planes, which include the paraboloid. For the paraboloid, the perimeter-minimizing region is a geodesic disc centered at the vertex of the paraboloid, where the Gaussian curvature is largest. The example of the paraboloid is also included in Theorem 1.2 in \cite{MHH}, where they study rotationally symmetric planes with its Gaussian curvature a strictly decreasing function from an origin. In this case, the boundaries of the perimeter-minimizing regions are again circles, although in general the region that these circles bound may be the complement of a geodesic disc or an annulus. Using a geometric measure theory and calculus of variations approach, Ritor\'e \cite{Ri} also proves an isoperimetric inequality for some complete surfaces, with rotational symmetry, including closed surfaces with equatorial symmetry, where the Gaussian curvature is monotonic away from the equator. Another general isoperimetric inequality for surfaces with variable curvature is given by Topping \cite{To}, and see the comprehensive surveys by Osserman \cite{Oss} and Howards, Hutchings, and Morgan \cite{HHM} for many other results and references. 

In this paper, we study surfaces which have constant positive Gaussian curvature, away from a finite number of singular, non-smooth points. These surfaces are formed from the double of a convex spherical polygon on the unit sphere $\mathbb{S}^2$. A special case is when the polygon is a spherical lune, in which case the surface has an explicit isometric embedding as a surface of revolution (see Definition \ref{defn:football} and Lemma \ref{lem:football} below). This special case then formally fits into the framework of \cite{MHH} and \cite{Ri}: The surface has rotational symmetry with the Gaussian curvature identically equal to $1$ except for two non-smooth tips where the curvature is undefined, but the tips can be thought of as providing an increase to the total Gaussian curvature. We use this framework to {give a new proof of} an isoperimetric inequality for these surfaces, with equality only attained when the surface is formed from spherical lunes and the region is a geodesic disc centered at a tip (see Theorem \ref{thm:iso}). {This inequality, and the case of equality, also follows directly by a straightforward approximation argument from the work of B\'erard, Besson, and Gallot in \cite{BBG} and Milman in \cite{EM1}, where they prove generalizations of the L\'evy-Gromov isoperimetric inequality \cite{LG1}. See the discussion after Theorem \ref{thm:iso} for more details.}

An important application of many isoperimetric inequalities is the relationship with Dirichlet eigenvalues via the Faber-Krahn theorem, \cite{Fa}, \cite{Kr}. This is carried out in a very general setting by Chavel, \cite{Ch}, where it is shown how combining an appropriate isoperimetric inequality with the co-area formula can provide a lower bound on the first Dirichlet eigenvalue of a region by that of a geodesic disc of the same area. Together with a careful analysis of the first Dirichlet eigenvalue of spherical caps, this has been for example used in \cite{BKP} to prove that two hemispheres provides the unique partition of the sphere that minimizes a certain function of the first Dirichlet eigenvalue of each piece. This function is related to harmonic, homogeneous, functions supported on cones generated by the subsets of the sphere, and is of importance in the regularity theory for minimizers of two-phase free boundary problems, \cite{ACF}. Since our surfaces are formed from two copies of a convex spherical polygon, a direct consequence of our isoperimetric inequality and a Faber-Krahn theorem is eigenvalue estimates for subsets of the polygon with mixed Dirichlet-Neumann boundary conditions (see Corollary \ref{cor:DN}).

\subsection*{Statement of Results}

As mentioned above, we will study an isoperimetric inequality and Faber-Krahn eigenvalue estimates for surfaces formed by copies of spherical lunes, and more generally (geodesically) convex spherical polygons. For each $a$, with $0<a \leq 1$, we set $\Omega_a$ to be a spherical lune on the unit sphere $\mathbb{S}^2$, with interior angles equal to $\pi a$. In particular, $\Omega_a$ has surface area $2a\pi$, and we define $\mathcal{P}_a$ to be the set of convex spherical polygons on $\mathbb{S}^2$ of the same surface area. For each such $P\in\mathcal{P}_a$, we will construct a closed (non-smooth) manifold $\tilde{P}$ consisting of gluing two copies of $P$ along its common boundary. For the spherical lune, this manifold will be constructed by first showing that $\Omega_a$ has an isometric embedding given by a surface with boundary contained in a plane and equal to a convex curve, which is smooth except possibly at two points. The closed manifold is then given by gluing two copies of this embedding across this plane, and this will lead to a smooth surface of revolution, except for two non-smooth tips. See Lemma \ref{lem:football} and its proof for the details of this construction. In particular, in this case of the spherical lune $\Omega_a$, we will see that an isometric embedding of the closed manifold is given by the following surface of revolution.
\begin{defn} \label{defn:football}
For each $0<a\leq 1$, define the surface of revolution $S_{a}$ by
\begin{align*}
S_{a} = \left\{(g(u),h(u)\cos(v),h(u)\sin(v)): -\tfrac{1}{2}\pi \leq u \leq \tfrac{1}{2}\pi, 0\leq v \leq 2\pi\right\}.
\end{align*}
Here the functions $g(u)$ and $h(u)$ are given by
\begin{align*}
g(u) = \int_{0}^{u}\left(1-a^2\sin^2(t)\right)^{1/2}\, dt, \qquad h(u) = a \cos(u).
\end{align*}
For each $b$, $-\tfrac{1}{2}\pi < b < \tfrac{1}{2}\pi$, we also define the cap $U_{a,b}$ given by
\begin{align*}
    U_{a,b} = \left\{(g(u),h(u)\cos(v),h(u)\sin(v)): b \leq u \leq \tfrac{1}{2}\pi, 0\leq v \leq 2\pi\right\}.
\end{align*}
\end{defn}
Note that when $a=1$, the defining function for $g(u)$ reduces to $g(u) = \sin(u)$, and so $S_1$ is the unit sphere $\mathbb{S}^2$. For $0<a<1$, $S_a$ is a smooth surface of revolution away from two singular points on the axis of rotation, and has the following properties.
\begin{lem} \label{lem:football}
The surfaces $S_a$ have the following properties:
\begin{enumerate}

  \item[i)]  Setting $a^* = g(\pi/2)$,  $S_{a}$ has constant Gaussian curvature $1$ everywhere, except possibly  at the two \textit{tips} at $(\pm a^*,0,0)$. The surfaces $S_{a}$ are the only simply connected surfaces of revolution with this constant Gaussian curvature property.
    
    \item[ii)] The metric on $S_{a}$ is given by $\,d u^2 + a^2\cos^2(u) \,d v^2.$

    \item[iii)]  The surface area of $S_a$ is equal to $4a\pi$.

    \item[iv)] The surface $S_a$ is an isometric embedding of the double of the lune $\Omega_a$, formed by gluing two copies of an isometric embedding of $\Omega_a$ along its boundary. 

\end{enumerate}

Moreover, setting $L$ to be the boundary length of a cap $U_{a,b}$, enclosing surface area $A$, it satisfies the equation
\begin{align*}
    L^2 = A(4a\pi - A).
\end{align*}
\end{lem}
 
\begin{figure}\begin{center}
\begin{tikzpicture}
\fill[blue!10](3,1.1) .. controls (2.5,1) and (2.7,3) .. (3,2.87);
\fill[blue!10] (3,1.14) .. controls (5.3,1.98) and (5.75,1.92) .. (3,2.85);
\draw [ultra thick, blue](3,1.1) .. controls (2.6,1) and (2.7,3) .. (3,2.87);
\draw [dashed, blue] (3,1.15) .. controls (3.1,1) and (3.4,3) .. (3,2.85);
\draw [ultra thick] (-5,2) .. controls (-1,4) and (1,4) .. (5,2);
\draw [ultra thick](-5,2) .. controls (-1,0) and (1,0) .. (5,2);
\draw [ultra thick](0,0.5) .. controls (-0.6,0.4) and (-0.6,3.7) .. (0,3.49);
\draw [dashed](0,0.52) .. controls (0.5,0.4) and (0.6,3.7) .. (0,3.47);
\node [right] at (5,2) {$(a^*,0,0)$};
\node [left] at (-5,2) {$(-a^*,0,0)$};
\node [below right, blue] at (2.8,1.1) {$u=b$};
\node [above] at (0,3.55) {$(0,a,0)$};
\node [below] at (0,0.45) {$(0,-a,0)$};
\node [above right, blue][font = \Large] at (3.5,2.8) {$U_{a,b}$};
\draw [thick, blue] (4,2.84) -- (3.6,2);

\end{tikzpicture}
\end{center}
\caption{The surface of revolution $S_a$ and cap $U_{a,b}$} \label{fig:football}
\end{figure}
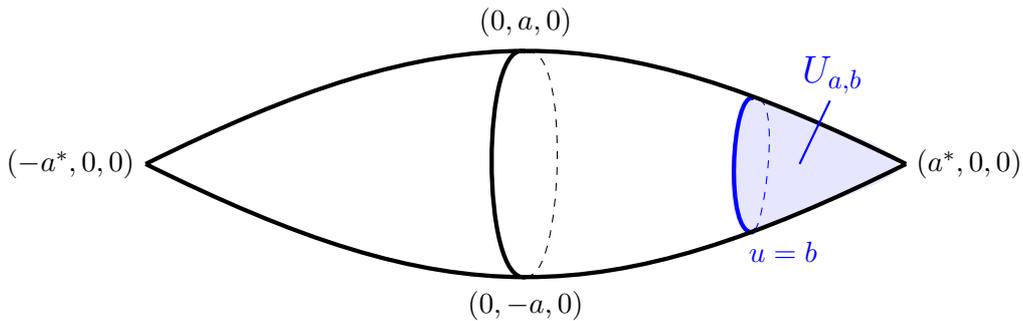

See Figure \ref{fig:football} for a visual representation of the surfaces $S_{a}$ and caps $U_{a,b}$.

The surfaces $S_a$ are given in Section 5.7 in \cite{ON}, where all surfaces of revolution of  constant positive Gaussian curvature are constructed. We will prove the remaining properties in Lemma \ref{lem:football} at the end of  Section \ref{sec:iso}. 
\begin{rem} \label{rem:tilde-P}
By Lemma \ref{lem:football} iv), for a general convex spherical polygon $P$, with $k$ vertices and interior angles $\theta_1,\ldots,\theta_k$, the resulting doubled manifold $\tilde{P}$ will have $k$ non-smooth points. In a neighborhood of such a point, corresponding to an interior angle $\theta_j$ on $P$, the manifold $\tilde{P}$ has an isometric embedding given by part of a surface $S_{a_j}$, with $a_j$ chosen so that the interior angle of the lune $\Omega_{a_j}$ is equal to $\theta_j$. 
\end{rem}
Informally, Theorem \ref{thm:iso} below states that for each $P\in \mathcal{P}_a$, the perimeter of any region on $\tilde{P}$ enclosing a region of surface area $A$, cannot be less than that of a cap $U_{a,b}$ on $S_a$ enclosing the same area. We will give a proof of this isoperimetric inequality for the following regions. 

\begin{defn} \label{defn:polygon}
Given $P\in\mathcal{P}_a$, the set of regions $\mathcal{U}_{\tilde{P}}$ on the doubled manifold  $\tilde{P}$ is defined as follows: A set $U$ is in $\mathcal{U}_{\tilde{P}}$ if its boundary consists of a finite number of disjoint, smooth, simple, closed curves $C_1,\ldots,C_m$ on $\tilde{P}$. That is, each curve $C_j$ is disjoint and a smooth embedding from $\mathbb{S}^1$ to $\tilde{P}$.
\end{defn}
We now state the isoperimetric inequality for the surfaces $S_a$ and $\tilde{P}$. 
\begin{thm}[Isoperimetric inequality] \label{thm:iso}
Given $a$ with $0<a\leq1$, a spherical convex polygon $P\in\mathcal{P}_a$, and a region $U\in \mathcal{U}_{\tilde{P}}$, let $L$ be the total perimeter of the boundary of $U$ and $A$ its surface area on $\tilde{P}$. Then, the isoperimetric inequality
\begin{align*}
    L^2 \geq A(4a\pi-A)
\end{align*}
holds, with equality if and only if $\tilde{P} = S_a$, and $U$ is a cap on $S_a$ enclosing surface area $A$.
\end{thm}
\begin{rem}
In the special case where $a=1$, then $P$ is a hemisphere and $\tilde{P} = \mathbb{S}^2$. In this case, Theorem \ref{thm:iso} reduces to the classical isoperimetric inequality for the sphere, stating that any perimeter-minimizing region must be a spherical cap. Therefore, from now on we will assume that $0<a<1$, and so in particular $S_{a}$ will be non-smooth at precisely two points.
\end{rem}
\begin{rem} \label{rem:quant1}
Let $\delta = \delta(P)$ be the difference between the smallest interior angle of $P$ and that of the lune $\Omega_{a}$ of the same area. Since $P$ is convex, note that $\delta(P)>0$ whenever $P\neq S_a$. Then, in the course of proving Theorem \ref{thm:iso}, we will show that, for regions $U\in\mathcal{U}_{\tilde{P}}$, the quantitative statement $L^2\geq (1+f(\delta))A(4a\pi-A)$ holds, for a function $f$ satisfying $f(\delta)>0$ for $\delta>0$. {Moreover, for sufficiently small $L$, depending on the distance between the vertices of the polygon $P$, the isoperimetric minimizers are caps centered around the tip corresponding to the smallest interior angle of $P$. In \emph{\cite{Mo1}}, an analogous property is shown for polytopes in any dimension, where, for small volume, geodesic balls about some vertex minimize perimeter.}
\end{rem}

{Theorem \ref{thm:iso} also gives an isoperimetric inequality for geodesically convex subsets $W$ of the sphere. In the corollary below $\mathcal{U}_{W}$ is equal to those regions on $W$ with boundary consisting of a finite number of disjoint, smooth, simple curves which are either closed or touch the boundary of $W$ at two points.
\begin{cor} \label{cor:approx}
Let $W$ be a geodesically convex subset of $\mathbb{S}^2$, with surface area $2a\pi$. For $U\in\mathcal{U}_{W}$, letting $L$ be the total perimeter of the boundary of $U$, and $A$ its surface area on $W$, the inequality
\begin{align*}
    L^2 \geq A(2\pi a - A)
\end{align*}
holds, with equality if and only if $W$ is the lune $\Omega_{a}$, and $U$ is the subset of the lune with boundary equal to a curve of constant latitude.
\end{cor}}

{By approximation, Theorem \ref{thm:iso} is also contained in  earlier work for smooth surfaces
: The L\'evy-Gromov inequality \cite{LG1} gives an isoperimetric inequality for compact, smooth manifolds with a lower bound on the Ricci curvature, by comparing to the model spaces $\mathbb{M}_{\kappa}$ of constant curvature $\kappa$. In particular, when $\kappa=1$, the model space is the unit sphere. In Theorem 2 in \cite{BBG}, B\'erard, Besson, and Gallot prove an improvement of this inequality when the diameter of the manifold is strictly less than that of the model space. In the setting of the sphere, this theorem ensures that for smooth manifolds, with Ricci curvature bounded below by $1$, and diameter strictly less than $\pi$, then the strict inequality from Theorem \ref{thm:iso} holds. In our setting, since $P$ is a convex polygon, the doubled surface $\tilde{P}$ has diameter strictly less than $\pi$ whenever $P$ is not a lune. Moreover, by a direct construction and calculation, the surface $\tilde{P}$ can be smoothly approximated by surfaces $\tilde{P}_{\eps}$ with curvature bounded below by $1$ (see Lemma \ref{lem:smoothed} and Remark \ref{rem:smoothed} below). Theorem 2 in \cite{BBG} therefore gives a strict isoperimetric inequality for each $\tilde{P}_{\eps}$, and by taking the limit, the strict inequality in Theorem \ref{thm:iso} is thus also contained in this theorem. To further generalize and improve the L\'evy-Gromov inequality and the inequality in \cite{BBG}, sharp isoperimetric inequalities have been proved for metric measure spaces with a Curvature-Dimension-Diameter condition providing a lower bound on a (generalized) Ricci curvature tensor, by Milman (see Corollary 1.4 in \cite{EM1}, and Section 6 in the same paper for an approximation procedure in the non-smooth case) and, in a non-smooth setting, by Cavalletti and Mondino, \cite{CM1}.

}

\subsection*{Application to eigenvalues on the polygon $P$ and its double $\tilde{P}$}

\begin{figure}   
    \centering
\begin{tikzpicture}

      \def\r{3}

    \draw circle (\r);
    \draw[dashed] ellipse (\r{} and \r/3);
      \fill[blue!10] (-0.83,1.63) .. controls (-0.2,1.38) and (0.2,1.38) .. (0.83,1.63);
    \fill[blue!10] (-0.83,1.63) .. controls (-0.001,3.36) and (0.001,3.36) .. (0.83,1.63);
    \draw[ultra thick] (0,-3) .. controls (-1.6,-1) and (-1.6,1) .. (0,3);
    \draw[ultra thick] (0,-3) .. controls (1.6,-1) and (1.6,1) .. (0,3);
    \draw[ultra thick, blue] (-0.85,1.6) .. controls (-0.2,1.4) and (0.2,1.4) .. (0.85,1.6);
    \draw[thick,blue] (0,2.3) -- (1.2,3.2);
    \node[right, blue][font = \Large] at (1.2,3.2) {$\Omega_{a,b}$};
    \node[left][font = \Large] at (-2.2,2.2) {$\mathbb{S}^2$};
    \node[right, blue] at (0.9,1.65) {$\phi = \tfrac{\pi}{2}+b$};

\end{tikzpicture}
    \caption{The spherical lune $\Omega_{a}$ and subset $\Omega_{a,b}$.} \label{fig:lune}
  
\end{figure}
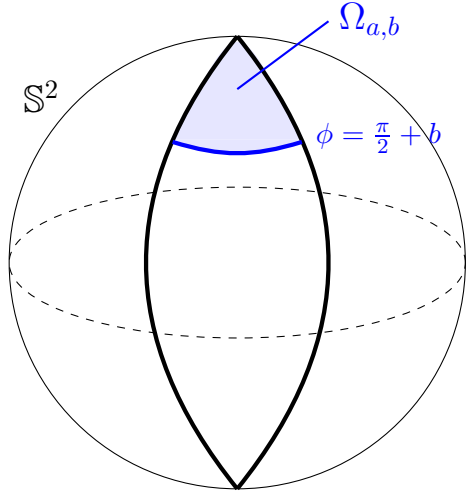

We can use Theorem \ref{thm:iso} to study the first eigenvalue of regions on the polygon $P$ and manifold $\tilde{P}$, with certain boundary conditions. Given $P\in\mathcal{P}_{a}$, and a region $U\in\mathcal{U}_{\tilde{P}}$, we define $\lambda(U)$ to be the first Dirichlet eigenvalue of $U$. That is,
\begin{align*}
    \lambda(U) = \inf\left\{\frac{\iint_{U}\left|\nabla_{g_{a}}w\right|^2\,d\sigma_{a} }{\iint_{U} w^2 \,d\sigma_{a}} : w\in C^{\infty}(U), w\big|_{\pa U} = 0 \right\}.
\end{align*}
Here $\nabla_{g_{a}}$ is the gradient on $\tilde{P}$, and $\,d\sigma_{a}$ the surface measure. In particular, when $\tilde{P}=S_{a}$ and $U = U_{a,b}$ is a cap, then by Lemma \ref{lem:football} ii), 
\begin{align} \label{eqn:eigenvalue1}
    \lambda(U_{a,b}) = \inf\left\{\frac{\int_{b}^{\pi/2}\int_{0}^{2\pi}\left[\left(\frac{\pa w}{\pa u}\right)^2+\frac{1}{a^2\cos^2(u)}\left(\frac{\pa w}{\pa v}\right)^2\right]a\cos(u)\,du\,dv }{\int_{b}^{\pi/2}\int_{0}^{2\pi} w^2\, a\cos(u)\,du\,dv} : w\in C^{\infty}(U_{a,b}), w(b,v) \equiv 0 \right\}.
\end{align}
Since $U_{a,b}$ is rotationally symmetric, and its first Dirichlet eigenvalue is simple, the corresponding eigenfunction must be independent of $v$. Therefore, from \eqref{eqn:eigenvalue1}, $\lambda(U_{a,b})$ is independent of $a$, and equals $\lambda(U_{1,b})$, the first Dirichlet eigenvalue of a spherical cap on $\mathbb{S}^2$. We now state a consequence of the isoperimetric inequality concerning these eigenvalues.
\begin{thm}[Faber-Krahn theorem] \label{thm:FK}
Given $a$ with $0<a\leq 1$, a spherical convex polygon $P\in\mathcal{P}_{a}$, and a region $U\in\mathcal{U}_{\tilde{P}}$, let $U_{a,b}$ be the cap on $S_a$ with the same surface area as $U$. Then, the first Dirichlet eigenvalue of $U$ satisfies
\begin{align*}
    \lambda(U) \geq \lambda(U_{a,b}) = \lambda(U_{1,b}),
\end{align*}
with equality if and only if $\tilde{P} = S_{a}$ and $U=U_{a,b}$ is a cap on $S_{a}$.
\end{thm}
\begin{rem} \label{rem:quant2}
As before, let $\delta = \delta(P)\geq0$ be the difference between the smallest interior angle of $P$ and that of the lune $\Omega_{a}$ of the same area. Then, the proof of Theorem \ref{thm:FK} will imply the quantitative statement $\lambda(U)\geq (1+g(\delta))\lambda(U_{a,b})$, for a function $g$ satisfying $g(\delta)>0$ for $\delta>0$.
\end{rem}
In Section \ref{sec:FK} we will explain how the proof of Theorem \ref{thm:FK} follows from the isoperimetric inequality in Theorem \ref{thm:iso} in an analogous way to how the classical Faber-Krahn theorem, \cite{Fa}, \cite{Kr}, follows from the classical isoperimetric inequality. For now, we explain a consequence of this theorem to eigenvalues of domains on geodesically convex subsets of $\mathbb{S}^2$: Let $W$ be a closed, geodesically convex set on $\mathbb{S}^2$, and let $V\subset W$ be a subset with smooth boundary. We define $\mu(V)$ to be the first eigenvalue of $V$, with Neumann boundary conditions on $\pa V\cap\pa W$ and Dirichlet boundary conditions on the rest of $\pa V$. One special case is when $W$ is a spherical lune $\Omega_{a}$, and $V = \Omega_{a,b}$ is the subset of the lune with boundary equal to a curve of constant latitude, $\phi = \tfrac{\pi}{2}+b$ (see Figure \ref{fig:lune}). In particular, $\Omega_a$ and $\Omega_{a,b}$ have isometric embeddings equal to one half of the surface $S_{a}$ and cap $U_{a,b}$ respectively. Via a doubling and approximation argument, in Section \ref{sec:FK} we will establish the following corollary of Theorem \ref{thm:FK}.

\begin{cor} \label{cor:DN}
Given a closed, geodesically convex set $W$ on $\mathbb{S}^2$ and a subset $V\subset W$ with a smooth boundary, let $\Omega_{a}$ and $\Omega_{a,b}$ be chosen to have the same area as $W$ and $V$ respectively. Then,
\begin{align*}
    \mu(V) \geq \mu(\Omega_{a,b}) = \mu(\Omega_{1,b}) = \lambda(U_{1,b}),
\end{align*}
with equality if and only if $W = \Omega_{a}$ and $V=\Omega_{a,b}$. 
\end{cor}
Given such a convex set $W$ and subset $V$, we use $\mu(V)$ to define the characteristic exponent of $V$,
\begin{align} \label{eqn:char1}
    \alpha(V) = -\tfrac{1}{2}+\sqrt{\tfrac{1}{4}+\mu(V)}.
\end{align}
This non-linear function of $\mu(V)$ is also the  homogeneity of a positively homogeneous, harmonic function on the cone generated by $V$, with Neumann boundary conditions on the part of the boundary of the cone generated by $\pa W$ and Dirichlet boundary conditions otherwise. Up to a scalar multiple, this harmonic function is unique. Beckner-Kenig-Pipher, \cite{BKP}, fully analyzed the eigenvalues $\lambda(U_{1,b})$ of the spherical caps $U_{1,b}$ in this context. They showed that for
\begin{align} \label{eqn:char2}
    \alpha(U_{1,b}) = -\tfrac{1}{2}+\sqrt{\tfrac{1}{4}+\lambda(U_{1,b})},
\end{align}
the inequality
\begin{align} \label{eqn:BKP}
    \alpha(U_{1,b}) + \alpha(U_{1,-b}) \geq 2
\end{align}
holds for $-\tfrac{\pi}{2} \leq b \leq \tfrac{\pi}{2}$. Moreover, they proved that equality holds if and only if $b=0$, so that the spherical caps $U_{1,b}$ and $U_{1,-b}$ are both hemispheres. Corollary \ref{cor:DN} therefore has the following consequence.
\begin{cor} \label{cor:FH}
Given a closed, geodesically convex set $W$ on $\mathbb{S}^2$ and a subset $V\subset W$ with a smooth boundary, the characteristic exponents of $V$ and $W\backslash V$ satisfy
\begin{align*}
    \alpha(V) + \alpha(W\backslash V)\geq 2.
\end{align*}
Moreover, we have equality only when $W$ is a lune $\Omega_{a}$, and $V = \Omega_{a,0}$.
\end{cor}
Corollary \ref{cor:FH} is a Dirichlet-Neumann version of the Friedland-Hayman inequality, \cite{FH}, and has been proved using different techniques coming from optimal transport and Caffarelli's contraction theorem (\cite{Ca1}, \cite{Ca2}) in \cite{BeJ}. The original Friedland-Hayman inequality and this version play an important role in the interior and boundary regularity of two-phase free boundary problems, see for example \cite{ACF}, \cite{BJR}.

\textbf{Acknowledgements.} {The authors would like to thank Frank Morgan and Emanuel Milman for their extremely useful comments and feedback, and in particular for explaining the results in \cite{EM1} and references therein.} The authors were supported by NSF Grant DMS-2042654.


\section{Proof of the isoperimetric inequality} \label{sec:iso}

In this section, we will give a proof of Theorem \ref{thm:iso}. We will split into two cases, by first proving the theorem for the surfaces $S_a$, and then extend to a general doubled polygon $\tilde{P}$. In both cases, key ingredients in the proof are the following isoperimetric inequalities. The first concerns a general inequality for simply-connected domains on surfaces with an upper bound on the Gaussian curvature.
\begin{thm}[Bol-Fiala inequality \cite{Bol}, \cite{Fia}, see Theorem 4.3 in \cite{Oss}] \label{thm:BF}
Let $D$ be a simply-connected region on a surface $S$, contained in a region where the Gaussian curvature $K$ satisfies $K\leq K_0$. Then, the perimeter $L$ and area $A$ of $D$ satisfy the inequality
\begin{align*}
    L^2\geq 4\pi A-K_0A^2.
\end{align*}
\end{thm}
Note that when $K = K_0 \equiv 1$, this is the isoperimetric inequality for the sphere. The other two results also concern isoperimetric inequalities for surfaces with prescribed information about their Gaussian curvature.
\begin{thm}[Morgan, Hutchings, Howards, \cite{MHH}, Theorem 2.1] \label{thm:Morgan1}
Let $S$ be a smooth Riemannian surface (either with or without boundary), with a smooth, rotationally invariant metric, such that the Gauss curvature is a non-increasing function of the distance from an origin. Then, among disjoint unions of embedded discs of a given area, the geodesic disc centered at the origin minimizes the perimeter. Moreover, unless the Gaussian curvature of a region containing the origin is constant, this is the unique minimizer.
\end{thm}
\begin{thm}[Morgan, Hutchings, Howards, \cite{MHH}, Theorem 2.4] \label{thm:Morgan2}
Let $S$ be a smooth Riemannian surface (either with or without boundary), and let $G(t)$ be the supremum of the total Gaussian curvature over regions of area $t$ on $S$. Let $L$ be the length of the boundary of a finite number of embedded discs, enclosing total area $A$. Then,
\begin{align*}
L^2 \geq 4\pi A - 2\int_{0}^{A}G(t) \, dt.
\end{align*} 
Moreover, one case of equality holds for a geodesic disc centered at the origin in a surface of revolution of non-increasing Gaussian curvature.
\end{thm}

\subsection{The proof for $S_a$}

We now prove Theorem \ref{thm:iso} for the lune $P = \Omega_a$, with corresponding surface $\tilde{P} = S_a$. Since the surfaces $S_a$ are non-smooth at the tips, when the region $U$ contains some of these tips, we will see that it is often more convenient to work with a smoothed version of $S_a$. However, in order to apply Theorems \ref{thm:Morgan1} and \ref{thm:Morgan2}, we want to maintain control on the Gaussian curvature of the smoothed surfaces.
\begin{lem} \label{lem:smoothed}
Given $\eps>0$, there exists a smooth surface $S_{a,\eps}$ with the following properties:
\begin{enumerate}

    \item[i)] $S_{a,\eps}$ is a surface of revolution
    \begin{align*}
        S_{a,\eps} = \left\{(g_{\eps}(u),h_{\eps}(u)\cos(v),h_{\eps}(u)\sin(v)): -\tfrac{1}{2}\pi \leq u \leq \tfrac{1}{2}\pi, 0\leq v \leq 2\pi\right\},
    \end{align*}
    with $g_\eps(u) = g(u)$, $h_{\eps}(u) = h(u)$ for $-\tfrac{1}{2}\pi+\eps \leq u \leq \tfrac{1}{2}\pi -\eps$.
    
    \item[ii)] There exists an absolute constant $C$, independent of $\eps>0$, such that for all $u$,
    \begin{align*}
        |g_{\eps}(u) - g(u)|+ |h_{\eps}(u) - h(u)| \leq C\eps.
    \end{align*}
    \item[iii)] The surface $S_{a,\eps}$ is symmetric about $u=0$, and for $u\geq0$, the Gaussian curvature is a non-increasing function of the distance from the tip $(g_{\eps}(\pi/2),0,0)$.
    
\end{enumerate}
\end{lem}
\begin{rem} \label{rem:smoothed}
For each convex spherical polygon $P$, we also define corresponding smooth surfaces $\tilde{P}_{\eps}$: As described in Remark \ref{rem:tilde-P}, in a neighborhood of each non-smooth point of $\tilde{P}$, the surface $\tilde{P}$ is equal to part of a surface $S_{a_j}$. Therefore, the smooth surfaces $\tilde{P}_{\eps}$ can be defined using the same construction as for the surfaces $S_{a,\eps}$ with an appropriate choice of $a$, in order to smooth each tip.
\end{rem}
We will prove Lemma \ref{lem:smoothed} with Lemma \ref{lem:football} at the end of this Section, and first use them to prove Theorem \ref{thm:iso}. To prove Theorem \ref{thm:iso} for $\tilde{P} = S_a$, we first assume that the region $U$ on $S_a$ is a single embedded disc with smooth boundary. To prove the inequality we split into four cases depending on the the number of tips contained in $U$.

\begin{enumerate}

    \item[1.] Suppose first that $U$ does not contain either tip of $S_a$. Then, $U$ is a simply connected region on $S_a$ contained in a region where the Gaussian curvature $K$ satisfies $K\equiv 1$. By the Bol-Fiala inequality, Theorem \ref{thm:BF}, the perimeter of $U$ satisfies
    \begin{align*}
        L^2 \geq A(4\pi - A) > A(4\pi a - A),
    \end{align*}
    since $a<1$, and so we obtain the desired strict inequality.
    
    \item[2.] If both tips are contained in the interior of $U$, then neither are contained in the complement $S_a\backslash U$. By Lemma \ref{lem:football} iii), this complement has area $4\pi a-A$, and the same perimeter as $U$, and so by the Bol-Fiala inequality, we obtain
        \begin{align*}
        L^2 \geq (4\pi a-A)(4\pi - (4\pi a - A)) = (4\pi a - A)(A + 4\pi(1-a)) > (4\pi a - A)A,
    \end{align*}
    again using $a<1$, and this gives the strict inequality.
    
    \item[3.] Now suppose that exactly one tip is contained in the interior of $U$, say $(g(\pi/2),0,0)$. Let $U_{a,b}$ be a cap of $S_a$, with $b$ chosen so that $U_{a,b}$ has the same area $A$. This cap has perimeter $L_{b}$ satisfying $L_{b}^2 = A(4a\pi - A)$, and so we need to show that $L\geq L_{b}$, with equality only when $U = U_{a,b}$. For $\eps>0$ sufficiently small, the boundaries of $U$ and $U_{a,b}$ are contained in the part of $S_a$ with $-\tfrac{1}{2}\pi +\eps\leq u \leq \tfrac{1}{2}\pi - \eps$. So, replacing $S_a$ by its smoothed version $S_{a,\eps}$ from Lemma \ref{lem:smoothed}, we obtain new sets $U_{\eps}$, $U_{a,b,\eps}$ of the same area $A_{\eps}$, and with perimeters still equal to $L$ and $L_{b}$ respectively. Note that $U_{a,b,\eps}$ is a geodesic disc centered at $(g_{\eps}(\pi/2),0,0)$ on $S_{a,\eps}$. Therefore, applying Theorem \ref{thm:Morgan1} to a rotationally symmetric subset of $S_{a,\eps}$ containing $U_{\eps}$, $U_{a,b,\eps}$ and only the one tip $(g_{\eps}(\pi/2),0,0)$, we obtain $L > L_{b}$, with equality only when $U_{\eps}$ (and hence also $U$) is a cap. Alternatively, we may instead directly apply Theorem 3.5 from \cite{Ri} in order to conclude that the isoperimetric domains of $S_{a,\eps}$ are caps. However, the above method of proof using the results from \cite{MHH} will be instructive when we come to consider the surfaces $\tilde{P}$ below.

    \item[4.] The final case is when the boundary of $U$ passes through at least one of the two tips. Given $\delta>0$, we can perturb $U$ to form a new embedded disc $U_{\delta}$, containing neither tip, with area $A_{\delta}$ and perimeter $L_{\delta}$ satisfying
    \begin{align*}
        |A-A_{\delta}|<\delta, \quad |L-L_{\delta}| < \delta. 
    \end{align*}
    Then, by Case 1., we have
    \begin{align*}
        L_{\delta}^2 \geq A_{\delta}(4\pi-A_{\delta}) = A_{\delta}(4\pi a - A_{\delta}) + 4\pi(1-a) A_{\delta}.
    \end{align*}
    Therefore, by choosing $\delta>0$ sufficiently small, depending on $1-a>0$, we can ensure that $L^2>A(4\pi a- A).$
    
\end{enumerate}

This completes the case when $U$ consists of a single embedded disc. To deal with the general case, we will use the following lemma.
\begin{lem} \label{lem:sum}
Given $a$, with $0<a<1$, and $m\geq2$, let $L_1,L_2,\ldots,L_m>0$, $A_1,A_2,\ldots,A_m>0$ be two sequences of positive numbers satisfying $L_j^2\geq A_j(4\pi a - A_j)$ for all $1\leq j \leq m$. Then, setting $L = L_1+\cdots +L_m$, $A = A_1+\cdots+A_m$, we have
\begin{align*}
    L^2 > A(4\pi a - A).
\end{align*}
\end{lem}
\begin{proof1}{Lemma \ref{lem:sum}}
Using the inequalities $L_j^2 \geq A_j(4\pi a - A_j)$, we obtain
\begin{align*}
    L^2 = \bigg(\sum_{j=1}^{m}L_j\bigg)^2 > \sum_{j=1}^{m}L_j^2 \geq \sum_{j=1}^{m}A_j(4\pi a - A_j) & = 4\pi a\sum_{j=1}^mA_j - \sum_{j=1}^mA_j^2 \\
    & > 4\pi aA - \bigg(\sum_{j=1}^{m}A_j\bigg)^2 = 4\pi a A - A^2,
\end{align*}
giving the required inequality.
\end{proof1}

Given $U\in \mathcal{U}_{S_a}$, suppose that $U$ consists of a finite number of disjoint embedded discs, of areas $A_1,\ldots,A_m$, with perimeters $L_1,\ldots,L_m$. Then, by applying the single embedded disc argument we have $L_j^2\geq A_j(4\pi a - A_j)$, for $1\leq j \leq m$, and so by Lemma \ref{lem:sum} we have the strict inequality whenever $m\geq 2$.

Finally, for a general $U\in \mathcal{U}_{S_a}$ let  $U_1,\ldots,U_M$ be its connected components, each of area $A_1,\ldots,A_M$ and perimeter $L_1,\ldots,L_M$. Each connected component $U_j$ is either a single embedded disc, or else $U_j$ is the complement of $m\geq2$ disjoint embedded discs. In the first case, we have already shown that $L_j^2\geq A_j(4\pi a -A_j)$, with equality only when $U_j$ is a cap. In the second case, let $A_{j,1},\ldots,A_{j,m}$ be the areas of these discs (so that $A_{j,1}+\cdots+A_{j,m} = 4\pi a - A_j$), with perimeters $L_{j,1},\ldots,L_{j,m}$ (so that $L_{j,1}+\cdots+ L_{j,m} = L$). Then, for $1\leq k \leq m$, we have 
\begin{align*}
    L_{j,k}^2\geq A_{j,k}(4\pi a- A_{j,k}).
\end{align*}
 Applying Lemma \ref{lem:sum} once to $L_{j,1},\ldots,L_{j,m}$ and $A_{j,1},\ldots,A_{j,m}$ implies that
\begin{align*}
    L_j^2 > (4\pi a - A_{j})A_{j}.
\end{align*}
Applying the lemma again to $L_1,\ldots,L_M$ and $A_1,\ldots,A_M$, therefore implies that $L^2 \geq A(4\pi a - A)$, and we have the strict inequality unless $M=1$, and the one connected component $U_1$ is a cap. This completes the proof of Theorem \ref{thm:iso} when the surface $U$ is given by $U = S_a$.

\subsection{The proof for $\tilde{P}$}

We now prove the isoperimetric inequality for the surface $\tilde{P}$ coming from a general convex spherical polygon $P$ of $3$ or more sides. Let $\theta_i$, $1\leq i \leq n$, be the interior angles of $P$, with corresponding exterior angles $\psi_i = \pi-\theta_i$. Note that since $P$ is convex, the exterior angles $\psi_i$ are strictly positive. For $P\in\mathcal{P}_a$, the surface area of $P$ is equal to $2\pi a$, and so the convexity of $P$ ensures that all of the interior angles $\theta_i$ satisfy
\begin{align} \label{eqn:interior}
    \theta_i > \pi a.
\end{align} 
This is because the interior angles of the lune $\Omega_a$ of the same area are $\pi a$. Moreover, by the Gauss-Bonnet theorem, the exterior angles of $P$ satisfy
\begin{align} \label{eqn:GB}
\sum_{i=1}^{n}\psi_i = 2\pi - \text{ Area}(P) = 2\pi(1-a).
\end{align}
As for $S_a$, rather than working with the non-smooth surface $\tilde{P}$, it will often be more convenient to work with the smooth surface $\tilde{P}_{\eps}$ (see Lemma \ref{lem:smoothed} and Remark \ref{rem:smoothed}). This surface has $n$ smoothed tips, corresponding to the interior angles $\theta_1,\ldots,\theta_n$.  In order to use Theorem \ref{thm:Morgan2}, we need control on the supremum of the total Gaussian curvature over regions of area $t$ in $\tilde{P}_{\eps}$.  The lemma below informally follows from the  Gaussian curvature at each tip of $\tilde{P}$ being a delta function of weight $2(\pi-\theta_{j})$. 
\begin{lem} \label{lem:G-eps}
Given $\eps>0$, let $V$ be a region on $\tilde{P}_{\eps}$ with smooth boundary, and containing the tips of $\tilde{P}_{\eps}$ corresponding to the interior angles $\theta_{i_1},\ldots,\theta_{i_k}$. Then, there exists an absolute constant $C>0$ such that the Gaussian curvature $K_{\eps}$ of $\tilde{P}_{\eps}$ satisfies
\begin{align*}
    \int_{V}K_{\eps} \leq \emph{Area}(V) + 2\sum_{j=1}^{k}\left(\pi-\theta_{i_j}\right) + C\eps.
\end{align*}
\end{lem}
\begin{proof1}{Lemma \ref{lem:G-eps}}
By breaking $V$ into $k$ pieces, it is sufficient to prove the lemma when $V$ contains one tip, corresponding to an interior angle $\theta_1$. Let $V_{\eps}$ be the part of $V\subset\tilde{P}_{\eps}$ contained in the $\eps$-neighborhood of this tip where $\tilde{P}_{\eps}$ differs from $\tilde{P}$. Then, since the curvature of $K_{\eps}$ equals $1$ on $V\backslash V_{\eps}$,
\begin{align} \label{eqn:G-eps2}
    \int_{V}K_{\eps} = \int_{V \setminus V_{\eps}}K_{\eps} + \int_{V_{\eps}}K_{\eps}  = \text{Area}(V\setminus V_{\eps}) +\int_{V_{\eps}}K_{\eps}\leq \text{Area}(V)+ \int_{V_{\eps}}K_{\eps}.
\end{align}
To bound the total curvature of $V_{\eps}$, we can work with the surface $S_{a_1,\eps}$ from Lemma \ref{lem:smoothed} corresponding to the smoothed double of the spherical lune with interior angle $\theta_1$.  Since $S_{a_1,\eps}$ is a smooth surface for all $\eps>0$, the Gauss-Bonnet theorem implies that
\begin{align*}
    \int_{S_{a_1,\eps}}K_{\eps} = 4\pi.
\end{align*}
Moreover, the surface area of $S_{a_1,\eps}$ satisfies
\begin{align} \label{eqn:surface-area}
    \left|\text{Area}(S_{a_1,\eps}) - \text{Area}(S_{a_1})\right| =   \left|\text{Area}(S_{a_1,\eps}) - 4\pi a_1\right| \leq C\eps.
\end{align}
Here and below  $C$ is an absolute constant (independent of $\eps>0$), which may change from line-to-line. The curvature $K_{\eps}$ of $S_{a_1,\eps}$ is equal to $1$ outside of the $\eps$-neighborhoods of each tip. Denoting these $\eps$-neighborhoods by $W_{\eps}$, we therefore have
\begin{align} \label{eqn:GB1}
  4\pi =  \int_{W_{\eps}}K_{\eps} + \int_{S_{a_1,\eps}\setminus W_{\eps}}K_{\eps} 
   =
    \int_{W_{\eps}}K_{\eps}+ \int_{S_{a_1,\eps}\setminus W_{\eps}}1 
=
    \int_{W_{\eps}}K_{\eps}+
    \text{Area}(S_{a_1,\eps}\setminus W_{\eps}). 
\end{align}
Moreover, there exists a constant $C$ such that
\begin{align*}
    \text{Area}(S_{a_1,\eps}\setminus W_{\eps}) \geq  
    \text{Area}(S_{a_1,\eps}) - C\eps.
\end{align*}
Thus, by \eqref{eqn:GB1}
\begin{align*}
    \int_{W_{\eps}}K_{\eps} & = 4\pi - \text{Area}(S_{a_1,\eps}\setminus W_{\eps}) \leq 4\pi- \text{Area}(S_{a_1,\eps}) + C\eps
    \leq 4\pi-4\pi a + C\eps =
    4(\pi - \theta_1) + C\eps.
\end{align*}
Since $W_{\eps}$ consists of $\eps$-neighborhoods of two tips, and $V_{\eps}$ is contained in one of these neighborhoods, this implies that
\begin{align*}
     \int_{V_{\eps}}K_{\eps} =   \frac{1}{2}\int_{W_{\eps}}K_{\eps} \leq  2(\pi - \theta_1) + C\eps.
\end{align*}
Using this inequality in \eqref{eqn:G-eps2} completes the proof of the lemma.
\end{proof1}
Before we prove Theorem \ref{thm:iso} for $\tilde{P}$, we need one more lemma.
\begin{lem} \label{lem:Gen_P}
Let $\theta_{i_1},\ldots,\theta_{i_n}$ be the interior angles of $P$ in a given order. Then, for any $m$, with $0\leq m \leq n$, at least one of the the following two inequalities holds:
\begin{align*}
   \frac{\left(\sum_{j=1}^{m}\theta_{i_j} \right) - (m-1) \pi}{\pi} \geq a
\quad\text{ or }\quad
    \frac{\left(\sum_{j=m+1}^{n}\theta_{i_j} \right) - (n - m-1) \pi}{\pi} \geq a. 
\end{align*}
\end{lem}
\begin{proof1}{Lemma \ref{lem:Gen_P}}
Adding the left hand sides of the two inequalities, we get
\begin{align}\label{eqn:Gen-P1}
    \frac{\left(\sum_{j=1}^{n}\theta_{i_j}\right) -  n\pi +2\pi}{\pi}=  \frac{-\left(\sum_{j=1}^{n}\psi_{i_j}\right) + 2\pi}{\pi}.
\end{align}
where $\psi_{i_j} = \pi-\theta_{i_j}$ are the corresponding exterior angles. By \eqref{eqn:GB}, the right hand side of \eqref{eqn:Gen-P1} equals $2a$, and so in particular at least one of the two inequalities in the statement of the lemma must hold.
\end{proof1}

We now proceed with the proof of Theorem \ref{thm:iso} for $\tilde{P}$. As for $S_{a}$, we start by considering the case where the region $U$ is a single embedded disc on $\tilde{P}$, with smooth boundary not passing through any tips. Replacing $\tilde{P}$ by $\tilde{P}_{\eps}$ for $\eps>0$ sufficiently small, we then obtain a new region $U_{\eps}$ on $\tilde{P}_{\eps}$ with the same perimeter of $L$, and area $A_{\eps}$, with
\begin{align} \label{eqn:A-eps}
    \left|A-A_{\eps}\right| \leq C\eps.
\end{align}
Here $C$ is an absolute constant, independent of $\eps$, which again may change from line-to-line. Let $\theta_{i_1},\ldots,\theta_{i_n}$ be the interior angles of $P$, so that the region $U$ contains $m$ tips of $\tilde{P}$, with $0\leq m \leq n$, and corresponding to interior angles $\theta_{i_1},\ldots,\theta_{i_m}$. Recalling that $\tilde{P}_{\eps}$ agrees with $\tilde{P}$ outside of an $\eps-$neighborhood of each tip, let $V$ be any region on $\tilde{P}_{\eps}$, which does not contain the $\eps-$neighborhoods of the tips corresponding to the angles $\theta_{i_{m+1}},\ldots, \theta_{i_n}$. Then, by Lemma \ref{lem:G-eps}, the total curvature of $V$ on $\tilde{P}_{\eps}$ is bounded from above by
\begin{align*}
    \text{Area}(V) + 2\sum_{j=1}^{m}(\pi-\theta_{i_j}) + C\eps.
\end{align*}
Therefore, we can apply Theorem \ref{thm:Morgan2} to the resulting surface (with boundary) where the $\eps$-neighborhood of tips corresponding to the angles $\theta_{i_{m+1}},\ldots,\theta_{i_n}$ have been removed, with
\begin{align} \label{eqn:G-eps1}
G_{\eps}(t) \leq t + 2\pi\sum_{j=1}^{m}\left(1-\theta_{i_j}/\pi\right) + C\eps.
\end{align}
 We obtain the lower bound
\begin{align} \nonumber
    L^2 \geq 4\pi A_{\eps} - 2 \int_{0}^{A_{\eps}} G_{\eps}(t) \,dt & \geq 4\pi A_{\eps} - 2 \int_{0}^{A_{\eps}} t + 2\pi\sum_{j=1}^{m}\left(1-\theta_{i_j}/\pi\right) + C\eps   \,dt \\ \nonumber
    &= 4\pi A_{\eps} - A_{\eps}^2 -4\pi \left(m - \frac{\sum_{j=1}^{m}\theta_{i_j}}{\pi}\right)A_{\eps} - 2A_{\eps}C\eps \\ \label{eqn:interior2}
    &= 4\pi \left(\frac{\sum_{j=1}^{m}\theta_{i_j}}{\pi} - (m-1) \right)A_{\eps} -A_{\eps}^2 - 2A_{\eps}C\eps. 
\end{align}
We can also apply the same argument to the complement $\tilde{P}_{\eps}/U_{\eps}$, which has area $A'_{\eps} =\text{ Area}(\tilde{P}_{\eps}) - A_{\eps}$, and obtain
\begin{align*}
    L^2 &\geq 4\pi \left(\frac{\sum_{j=m+1}^{n}\theta_{i_j}}{\pi} - (n-m-1) \right)A'_{\eps} -(A'_{\eps})^2 - 2A'_{\eps}C'\eps.
\end{align*}
for another absolute constant $C'$. Lemma \ref{lem:Gen_P} thus implies that either
\begin{align*}
  L^2 \geq  4\pi \left(\frac{\sum_{j=1}^{m}\theta_{i_j}}{\pi} - (m-1) \right)A_{\eps} -A_{\eps}^2 - 2A_{\eps}C\eps &\geq 4\pi a A_{\eps} - A_{\eps}^2 - 2A_{\eps}C\eps.
\end{align*}
or 
\begin{align*}
  L^2\geq   4\pi \left(\frac{\sum_{j=m+1}^{n}\theta_{i_j}}{\pi} - (n-m-1) \right){A}'_{\eps} -({A}'_{\eps})^2 - 2{A}'_{\eps}C'\eps
    \geq  4\pi a {A}'_{\eps} - ({A}'_{\eps})^2 - 2{A}'_{\eps}C'\eps
\end{align*}
holds. Letting $\eps$ tend to $0$ and using \eqref{eqn:A-eps}, in either case we get the desired lower bound of
\begin{align*}
    L^2 &\geq 4\pi a A - A^2.
\end{align*}
\begin{rem} \label{rem:single}
If $U$ (or its complement) only contains one tip. corresponding to an angle $\theta$, then from \eqref{eqn:interior} and \eqref{eqn:interior2}, we immediately get the strict inequality $L^2>A(4\pi a - A)$.
\end{rem}
This completes the proof of the inequality when $U$ is a single embedded disc, with boundary not passing through any tips. When the boundary passes through some of the tips, we can proceed as we did in Case 4) for $S_a$, by perturbing $U$ in such a way that the sum of the  angles in either the new set or its complement strictly satisfy one of the inequalities in Lemma \ref{lem:Gen_P}. Finally, for a general region $U\in \mathcal{U}_{\tilde{P}}$, we can proceed as for $S_a$ using Lemma \ref{lem:sum} in order to establish the inequality.
\\
\\
We are left to prove that when $\tilde{P}\neq S_a$, we always have the strict inequality $L^2>A(4\pi a - A)$. From Remark \ref{rem:single}, when $\tilde{P}\neq S_a$ and $U$ only contains one tip, then we do get the strict inequality. In fact, for sufficiently small area $A$ (depending on the distance between the vertices of the polygon $P$), the minimizer can only contain at most one tip. By our isoperimetric inequality on $S_a$, the minimizer will therefore be a cap centered around the tip corresponding to the smallest interior angle, say $\theta_1$. Denoting $L(t)$ to be the shortest perimeter of regions in $\mathcal{U}_{\tilde{P}_{\eps}}$ enclosing area $t$, then for sufficiently small $t>0$, this ensures that
\begin{align} \nonumber
    L(t)^2 & \geq t\left(4\pi\tfrac{\theta_1}{\pi} -t\right) -C\eps > t\left(4\pi a-t\right), \\ \label{eqn:equality1} L(t)L'(t) & \geq 2\theta_1 - t-C\eps  = 2\pi -\left(t+2\pi\left(1-\tfrac{\theta_1}{\pi}\right)\right)-C\eps > 2\pi -\left(t+2\pi\left(1-a\right)\right).
\end{align}
Theorems \ref{thm:Morgan1} and \ref{thm:Morgan2} in \cite{MHH} are proved by integrating the inequality
\begin{align*}
    L(t)L'_{L^*}(t) \geq 2\pi - G(t),
\end{align*}
from $t=0$ to $t=A$. Here $L'_{L^*}(t)$ is the lower-left derivative of $L(t)$ (see (2) in \cite{MHH}, page 4894) and $G(t)$ is the supremum of the total Gaussian curvature of a region of area $t$. Therefore, instead using \eqref{eqn:equality1} for small values of $t$, and the upper bound on $G(t)$ from \eqref{eqn:G-eps1} otherwise, we get the strict inequality
\begin{align*}
    L^2 > A_{\eps}(4\pi a- A_{\eps})- 2A_{\eps}C\eps
\end{align*}
where the difference between the two sides is bounded below by a constant independent of $\eps>0$. Letting $\eps$ tend to $0$ therefore gives the strict inequality $L^2>A(4\pi a-A)$, and completes the proof of Theorem \ref{thm:iso}. Note that using \eqref{eqn:equality1} in this way also ensures that the quantitative inequality given in Remark \ref{rem:quant1} holds.

{\begin{proof1}{Corollary \ref{cor:approx}}
Given a geodesically convex set $W$ on $\mathbb{S}^2$ and $U\in \mathcal{U}_{W}$, if $W = \Omega_{a}$ is a lune, then the corollary follows immediately from Theorem \ref{thm:iso} by using $U$ to form a region on the doubled closed surface $\tilde{W} = S_a$ of twice the perimeter and area. 

If $W$ is not a lune, we first form a sequence of spherical polygons $P_k\subset W$, and corresponding closed surfaces $\tilde{P}_k$, such that $P_k$ converges to $W$ in Hausdorff measure. This leads to a sequence of regions $U_k$ on $P_k$, with perimeters and surface area converging to that of $U$. Since $W$ is convex but not a lune, it cannot contain antipodal points on $\mathbb{S}^2$, and thus the same is true for $P_k$. This gives a uniform lower bound of $\delta>0$ on $\delta_k = \delta(P_k)$, the difference between the smallest interior angle of $P_k$, and the lune $\Omega_{a_k}$ of the same area. Applying the quantitative isoperimetric inequality in Remark \ref{rem:quant1} to each closed surface $\tilde{P}_k$ and double of $U_k$, and then taking the limit as $k$ tends to infinity therefore proves the Corollary.
\end{proof1}
}

We are left to prove the properties of $S_{a}$ and its smoothed version given in Lemmas \ref{lem:football} and \ref{lem:smoothed}.

\begin{proof1}{Lemma \ref{lem:football}}
In Example 7.5 in Section 5.7 of \cite{ON}, all of the surfaces of revolution with constant positive curvature are constructed. These are given by
\begin{align*}
    \textbf{r}_{c}(u,v) = \left(g_{c}(u),h_{c}(u)\cos(v),h_{c}(u)\sin(v)\right),
\end{align*}
with 
\begin{align*}
    g_{c}(u) = \int_{0}^{u}\left(1-\tfrac{a^2}{c^2}\sin^2(t/c)\right)^{1/2}\,dt, \qquad h_{c}(u) = a\cos(u/c).
\end{align*}
Here $a$, $c >0$ are constants, and for all such $a$ and $c$, the surfaces have constant curvature $K = \tfrac{1}{c^2}$. Moreover, as shown in \cite{ON}, these surfaces are only simply connected for $0<a\leq c$, and for this range of $a$, the variables $u$ and $v$ are defined for $-\tfrac{1}{2}c\pi \leq u \leq \tfrac{1}{2}c\pi$ and $0 \leq v \leq 2\pi$. Therefore, setting $c = 1$, and $g(u) = g_1(u)$, $h(u) = h_1(u)$, this ensures that the surfaces $S_a$ from Definition \ref{defn:football} have constant Gaussian curvature equal to $1$ everywhere except at $(\pm g(\pi/2),0,0)$, and are the only simply connected surfaces of revolution with this property. This proves part i) of the lemma.
\\
\\
Setting $c=1$ and writing $\textbf{r}(u,v) = (g(u),h(u)\cos(v),h(u)\sin(v))$,  a direction calculation gives the expressions
\begin{align*}
      \mathbf{r}_{u}  & = \left( \sqrt{1-a^2\sin^2(u)}, -a\sin(u)\cos(v), -a\sin(u)\sin(v)\right), \\
    \mathbf{r}_{v} & = \left(0, -a\cos(u)\sin(v), a\cos(u)\cos(v)\right), \\
|\mathbf{r}_{u} \times \mathbf{r}_{v}| & = a\cos(u).
\end{align*}
This means that the surface $S_a$ indeed has the metric $du^2+a^2\cos^2(u)\,dv^2$, giving part ii) of the lemma.  
\\
\\
Moreover, the surface area of $S_{a}$ is equal to
\begin{align*}
    \int_{0}^{2\pi}\int_{-\pi/2}^{\pi/2}|\mathbf{r}_{u}\times\mathbf{r}_{v}|\, du\,dv  = \int_{0}^{2\pi}\int_{-\pi/2}^{\pi/2} a\cos(u) \, du\, dv = 4a\pi,
\end{align*}
proving iii).

The lune $\Omega_{a}$ has a parameterization given by
\begin{align*}
    \Omega_{a} = \left\{(\sin(\tilde{u}),\cos(\tilde{u})\cos(\tilde{v}),\cos(\tilde{u})\sin(\tilde{v}))\,:\, -\tfrac{1}{2}\pi \leq \tilde{u}\leq\tfrac{1}{2}\pi, 0 \leq \tilde{v} \leq a\pi \right\},
\end{align*}
with metric $d\tilde{u}^2 + \cos^2(\tilde{u})\,d\tilde{v}^2$. Therefore, using ii), the mapping $(u,v)\mapsto(\tilde{u},\tilde{v}) = (u,av)$ provides an isometry from one half of $S_{a}$, given by,
\begin{align*}
    \left\{(g(u),h(u)\cos(v),h(u)\sin(v)): -\tfrac{1}{2}\pi \leq u \leq \tfrac{1}{2}\pi, 0\leq v \leq \pi\right\},
\end{align*}
to the lune $\Omega_{a}$. In particular, the surface $S_{a}$ is indeed an isometric embedding of the closed manifold formed by gluing two copies of $\Omega_{a}$, giving iv) in the lemma.

Finally, suppose that a cap $U_{a,b}$ on $S_{a}$ encloses area $A$ and has boundary length $L$. By calculating the surface area of $U_{a,b}$, we see that $b$ and $A$ satisfy the equation
\begin{align*}
   A = \int_{0}^{2\pi}\int_{b}^{\pi/2}\ a\cos(u) \, du\,dv = 2a\pi(1-\sin(b)).
\end{align*}
Since every vertical cross-section (with $u$ fixed) of $S_{a}$ is a circle with radius $2a\pi\cos(u)$, the boundary length $L$ of the cap $U_{a,b}$ is therefore given by
\begin{align*}
    L^2 = (2a\pi\cos(b))^2 =  4a^2\pi^2(1-\sin^2(b)) = 4a^2\pi^2\left(1-\left(1-\tfrac{A}{2a\pi}\right)^2 \right),
\end{align*}
which simplifies to $L^2 = 4a\pi A-A^2$, and completes the proof of Lemma \ref{lem:football}.
\end{proof1}

We are left to prove Lemma \ref{lem:smoothed}. We prove this lemma by smoothing out the left and right tips of $S_a$ carefully in order to ensure that the Gaussian curvature is still an increasing function of the distance from the (smoothed) tips. 
\begin{proof1}{Lemma \ref{lem:smoothed}}
By translating and reparameterizing the surface $S_a$, we can write the left half of it as
\begin{align*}
    \textbf{r}(u,v) = \left(w(u), u\cos(v),u\sin(v)\right),
\end{align*}
for $0\leq v \leq 2\pi$, $0\leq u \leq a $, and a function $w(u)$ satisfying $0 = w(0) \leq w(u) \leq g(\pi/2)$. For $0<a<1$, the right derivative of $w(u)$ at $u=0$ is non-zero, reflecting the non-smoothness of $S_a$ at its two tips.
\\
\\
To obtain a smooth approximation of this part of $S_a$, we fix $\eps>0$ and consider the surface
\begin{align*}
    \textbf{r}_{\eps}(u,v) = \left(w_{\eps}(u), u\cos(v),u\sin(v)\right),
\end{align*}
for $0 \leq v \leq 2\pi$, $0 \leq u \leq a$, and a $C^2$-smooth function $w_{\eps}(u)$ with these two properties:
\begin{enumerate}
    \item[i)] $w_{\eps}(u) = w(u)$ for $\eps \leq u \leq a$;
    \item[ii)] $w_{\eps}(u) = b_0 + b_1u^2 + b_2u^4$ for $0 \leq u\leq \eps$, and for coefficients $b_0$, $b_1$, and $b_2$.
\end{enumerate}
The coefficients $b_0$, $b_1$ and $b_2$ are chosen so that
\begin{align} \label{eqn:smoothed0}
    w(\eps)=w_{\eps}(\eps), \quad  w'(\eps)=w'_{\eps}(\eps), \quad w''(\eps)=w''_{\eps}(\eps),
\end{align}
to ensure that $w_{\eps}$ is a $C^2$-smooth function. Solving these equations for $b_0$, $b_1$, and $b_2$ gives
\begin{align} \label{eqn:smoothed1}
    b_0 = \frac{8 w(\eps) - 5w'(\eps)\eps + w''(\eps)\eps^2}{8}, \quad b_1 = \frac{3w'(\eps) - w''(\eps)\eps}{4\eps}, \quad b_2 = \frac{-w'(\eps) + w''(\eps)\eps}{8\eps^3}.
\end{align}
Note that by the properties of the surfaces $S_a$, for $\eps>0$ sufficiently small, $w(\eps)$ and $w''(\eps)$ are both positive and bounded by a constant multiple of $\eps$, while $w'(\eps)$ is bounded above and below by a positive constant.  In particular, for $0 \leq u\leq \eps$, we therefore have $|w_{\eps}(u)|\leq C\eps$ for an absolute constant $C$. Defining the symmetric smoothed surface $S_{a,\eps}$ to be formed by two copies of $\textbf{r}_{\eps}(u,v)$, this ensures that i) and ii) in Lemma \ref{lem:smoothed} hold. These formulas also ensure that $b_1=b_1(\eps)>0$ and $b_2=b_2(\eps)<0$ for $\eps>0$ sufficiently small. 
\\
\\
To complete the proof of the lemma, we need to show that, for $v$ fixed, the Gaussian curvature of the surface $\textbf{r}_{\eps}(u,v)$ is a non-increasing function of $u$, with $0 \leq u \leq a$. We know that the curvature is identically equal to $1$ for $\eps \leq u\leq\eps$ by the properties of $S_{a}$, and so we only need to consider $0\leq u \leq a$. For this range of $u$, the Gaussian curvature is given by
\begin{align} \label{eqn:smoothed2}
    K(u) = \frac{u^{-1}w_{\eps}'(u)w_{\eps}''(u) }{(1+|w_{\eps}'(u)|^2)^2} = \frac{(2b_1+4b_2u^2)(2b_1+12b_2u^2)}{(1+(2b_1u+4b_2u^3)^2)^2},
\end{align}
and from \eqref{eqn:smoothed0} we know that $K(\eps) = 1$, matching the constant Gaussian curvature of $S_a$ away from the tips. By \eqref{eqn:smoothed1}, there exists a constant $C$, independent of $\eps>0$, such that 
\begin{align} \label{eqn:smoothed3}
    \left|b_1 - \tfrac{3}{4}\eps^{-1}w'(\eps) \right| \leq C\eps, \qquad  \left|b_2 + \tfrac{1}{8}\eps^{-3}w'(\eps) \right| \leq C\eps^{-1}.
\end{align}
These estimates guarantee that, for $\eps>0$ sufficiently small, the numerator in \eqref{eqn:smoothed2} is strictly decreasing for $0\leq u \leq \eps$, and since $K(\eps)=1\geq0$, this numerator must therefore also be positive for this range of $u$. 
The estimates in \eqref{eqn:smoothed3} also ensure that, for $\eps>0$ sufficiently small,
\begin{align} \label{eqn:smoothed4}
    w_{\eps}'(u) > 0, \quad w_{\eps}'''(u) <0, \quad w_{\eps}''(u)\geq w_{\eps}''(\eps) = w''(\eps) >0
\end{align}
for $0\leq u \leq \eps$. Therefore, the denominator of $K(u)$ in \eqref{eqn:smoothed2} is a strictly increasing, positive function. Putting everything together, we thus have that $K(u)$ is indeed a positive, decreasing function for $0\leq u \leq \eps$,  with $K(\eps) = 1$.
\end{proof1}

\section{Proof of the Faber-Krahn theorem and Corollaries \ref{cor:DN} and \ref{cor:FH}} \label{sec:FK}

In $2$ dimensions, the classical Faber-Krahn theorem states the following (see, for example, \cite{Ch}, Theorem 2, page 87).
\begin{thm}[Faber \cite{Fa}, Krahn \cite{Kr}] \label{thm:FK1}
Let $\mathbb{M}_{\kappa}$ be the complete, simply connected $n$-dimensional space of constant sectional curvature $\kappa$, and let $M$ be a complete $2$-dimensional Riemannian surface. For each open set $U$ consisting of a finite union of disjoint regular domains on $M$, let $D$ be the geodesic disc in $\mathbb{M}_{\kappa}$ of the same surface area. If for all such $U$, the boundary length of $U$ is greater than or equal to that of $D$, with equality if and only if $U$ and $D$ are isometric, then
\begin{align*}
\lambda(U)\geq\lambda(D).
\end{align*}
Here $\lambda(\cdot)$ are the respective first Dirichlet eigenvalues, and equality holds if and only if $U$ and $D$ are isometric.
\end{thm}
To prove Theorem \ref{thm:FK}, we follow the proof presented in \cite{Ch}, with the surface $\mathbb{M}_{\kappa}$ replaced by the surface $S_{a}$, chosen so that $S_a$ has the same surface area as the given doubled polygon $\tilde{P}$: The geodesic discs of $S_{a}$, centered at a tip $(a^*,0,0)$, are then precisely the caps $U_{a,b}$. Therefore, the isoperimetric inequality in Theorem \ref{thm:iso} ensures that the hypothesis on the boundary length of $U$ is satisfied. Theorem \ref{thm:FK} then follows by exactly replicating the proof of Theorem 2 on page 87 of \cite{Ch} by building a comparison function on $S_{a}$, with superlevel sets given by the caps $U_{a,b}$, and the use of the co-area formula. The quantitative eigenvalue lower bound in Remark \ref{rem:quant2} also follows by using the quantitative isoperimetric inequality in Remark \ref{rem:quant1} in the same argument.
\\
\\
We end by using Theorem \ref{thm:FK} to prove Corollaries \ref{cor:DN} and \ref{cor:FH}.
\begin{proof1}{Corollaries \ref{cor:DN} and \ref{cor:FH}}
We start by proving Corollary \ref{cor:DN}, and first consider the case when the geodesically convex set $W$ is a spherical polygon $P$. We let $U$ be the double of $V$ on the closed surface $\tilde{P}$. Then, by reflecting across the boundary of $W$, using the Neumann boundary conditions on $\pa W\cap\pa V$, the Dirichlet-Neumann eigenvalue $\mu(V)$ is a Dirichlet eigenvalue of $U$ on $\tilde{P}$. In particular, $\mu(V)\geq \lambda(U)$. Since the first Dirichlet eigenvalue of the cap $U_{a,b}$ is simple, it is rotationally symmetric, and so when $W$ is equal to a lune $\Omega_{a}$, and $V = \Omega_{a,b}$, this process must lead to the first Dirichlet eigenvalue of $U_{a,b}$. Therefore,
\begin{align*}
    \mu(\Omega_{a,b}) = \lambda(U_{a,b}) = \lambda(U_{1,b}).
\end{align*}
The inequality in the corollary then follows from Theorem \ref{thm:FK}. Moreover, if $\delta = \delta(P)\geq 0$ is the difference between the smallest interior angle of $P$ and that of $\Omega_a$, then 
\begin{align} \label{eqn:quant1}
    \mu(V) \geq (1+g(\delta))\mu(\Omega_{a,b}),
\end{align}
with $g(\delta)$ the function from Remark \ref{rem:quant2}.
\\
\\
For a general geodesically convex set $W$, we approximate the first Dirichlet-Neumann eigenvalue of a subset $V\subset W$ using this lemma:
\begin{lem} \label{lem:approx}
Let $V_k\subset V$ be subsets of $\mathbb{S}^2$, with boundaries of the following form: There exist sets $\pa V_k^N$, $\pa V_k^D$, $\pa V^N$, and $\pa V^D$, such that 
\begin{align*}
    \pa V_k = \pa V_k^N\cup\pa V_k^D, \qquad  \pa V = \pa V^N\cup\pa V^D,
\end{align*}
$\pa V_k^D\subset \pa V^D$, and $\pa V^N$ is Lipschitz. Then, denoting $\mu(V)$ to be the first eigenvalue of $V$, with Dirichlet boundary conditions on $\pa V^D$ and Neumann boundary conditions on $\pa V^N$, and likewise for the eigenvalue $\mu(V_k)$, we have
\begin{align*}
    \mu(V_k) \leq \mu(V)\left(1 + C\emph{Area}_{\mathbb{S}^2}(V\backslash V_k)\right).
\end{align*}
The constant $C>0$ depends only on the eigenvalue  $\mu(V)$ and the Lipschitz constant of $\pa V^N$.
\end{lem}
\begin{proof1}{Lemma \ref{lem:approx}}
Denote $u$ to be the corresponding first $L^2(V)$-normalized eigenfunction of $V$, and $\tilde{u}$ to be the reflection of $u$ across the Neumann boundary $\pa V^N$. Due to the Neumann boundary conditions on $\pa V^{N}$, and since $\pa V^N$ is Lipschitz, $\tilde{u}$ is then a solution of a non-degenerate elliptic equation in divergence form, with bounded, measurable coefficients on the double of $V$. These coefficients can be bounded in terms of $\mu(V)$ and the Lipschitz constant of $\pa V^N$, and so elliptic estimates (see, for example, Theorem 8.25 in \cite{GT}) ensure that $\tilde{u}$ and $u$ are bounded in $L^{\infty}$.

To obtain the upper bound on $\mu(V_k)$ given in the statement of the lemma, we use the variational formulation of the first eigenvalue. The restriction of $u$ to $V_k$, which we call $u_{k}$, is an admissible test function in this formulation, due to the structure of the boundaries and boundary conditions of $V_k$ and $V$. Therefore, letting $g$ be the round metric on $\mathbb{S}^2$, with surface measure $\,d\sigma_{g}$, 
\begin{align*}
    \mu(V_k) \leq \frac{\iint_{V_k} \left|\nabla_{g}u_k \right|^2\,d\sigma_{g}}{\iint_{V_k} \left|u_k \right|^2\,d\sigma_{g}} \leq \frac{\iint_{V} \left|\nabla_{g}u\right|^2\,d\sigma_{g}}{\iint_{V_k} \left|u \right|^2\,d\sigma_{g}} = \frac{\mu(V)}{\iint_{V_k}\left|u \right|^2\,d\sigma_{g}}.
\end{align*}
Using
\begin{align*}
    \iint_{V_{k}}\left|u \right|^2\,d\sigma_{g} = \iint_{V}\left|u \right|^2\,d\sigma_{g}  -  \iint_{V\backslash V_{k}}\left|u \right|^2\,d\sigma_{g} = 1- ||{u}||_{L^{\infty}(V)}\text{Area}(V\backslash V_{k}),
\end{align*}
and the boundedness of $u$, then gives the desired inequality.
\end{proof1}

To complete the proof of Corollary \ref{cor:DN}, let $W_k\subset W$ be a sequence of convex polygons on $\mathbb{S}^2$, converging to $W$ in Hausdorff measure, and let $V_k = V\cap W_k\subset V$. Since $V_k$ is a subset of a convex spherical polygon, we have $\mu(V_k) \geq \mu(\Omega_{a,b_k})$ where $\Omega_{a,b_k}$ has the same surface area as $V_k$. Applying Lemma \ref{lem:approx}, and letting $k$ tend to infinity, we therefore have the required inequality $\mu(V)\geq \mu(\Omega_{a,b})$. Moreover, if $W$ is not a lune $\Omega_{a}$, then we can choose the approximating sequence $W_k$, so that the difference between the smallest interior angle of $W_k$ and that of the lune $\Omega_{a_k}$ of the same area is bounded below by a positive constant, independently of $k$. The strict inequality in \eqref{eqn:quant1} therefore ensures that $\mu(V)>\mu(\Omega_{a,b})$ in this case.

This finishes the proof of Corollary \ref{cor:DN}. With $\alpha(V)$ and $\alpha(U_{1,b})$ as in \eqref{eqn:char1} and \eqref{eqn:char2}, Corollary \ref{cor:DN} ensures that
\begin{align*}
    \alpha(V) \geq \alpha(U_{1,b}),
\end{align*}
with equality only when $W=\Omega_{a}$ and $V = \Omega_{a,b}$. Corollary \ref{cor:FH} therefore follows immediately from \eqref{eqn:BKP}. 
\end{proof1}

\end{document}